\documentclass[10pt,draft]{amsart}
\usepackage{amsmath,amssymb,amsthm}
\begin{document}
\newtheorem{lem}{Lemma}[section]
\newtheorem{prop}{Proposition}[section]
\newtheorem{cor}{Corollary}[section]
\numberwithin{equation}{section}
\newtheorem{thm}{Theorem}[section]
\theoremstyle{remark}
\newtheorem{example}{Example}[section]
\newtheorem*{ack}{Acknowledgment}
\theoremstyle{definition}
\newtheorem{definition}{Definition}[section]
\theoremstyle{remark}
\newtheorem*{notation}{Notation}
\theoremstyle{remark}
\newtheorem{remark}{Remark}[section]
\newenvironment{Abstract}
{\begin{center}\textbf{\footnotesize{Abstract}}%
\end{center} \begin{quote}\begin{footnotesize}}
{\end{footnotesize}\end{quote}\bigskip}
\newenvironment{nome}
{\begin{center}\textbf{{}}%
\end{center} \begin{quote}\end{quote}\bigskip}

\newcommand{\triple}[1]{{|\!|\!|#1|\!|\!|}}
\newcommand{\xx}{\langle x\rangle}
\newcommand{\ep}{\varepsilon}
\newcommand{\al}{\alpha}
\newcommand{\be}{\beta}
\newcommand{\de}{\partial}
\newcommand{\la}{\lambda}
\newcommand{\La}{\Lambda}
\newcommand{\ga}{\gamma}
\newcommand{\del}{\delta}
\newcommand{\Del}{\Delta}
\newcommand{\sig}{\sigma}
\newcommand{\ome}{\omega}
\newcommand{\Ome}{\Omega}
\newcommand{\C}{{\mathbb C}}
\newcommand{\N}{{\mathbb N}}
\newcommand{\Z}{{\mathbb Z}}
\newcommand{\R}{{\mathbb R}}
\newcommand{\Rn}{{\mathbb R}^{n}}
\newcommand{\Rnu}{{\mathbb R}^{n+1}_{+}}
\newcommand{\Cn}{{\mathbb C}^{n}}
\newcommand{\spt}{\,\mathrm{supp}\,}
\newcommand{\Lin}{\mathcal{L}}
\newcommand{\SSS}{\mathcal{S}}
\newcommand{\F}{\mathcal{F}}
\newcommand{\xxi}{\langle\xi\rangle}
\newcommand{\eei}{\langle\eta\rangle}
\newcommand{\xei}{\langle\xi-\eta\rangle}
\newcommand{\yy}{\langle y\rangle}
\newcommand{\dint}{\int\!\!\int}
\newcommand{\hatp}{\widehat\psi}
\renewcommand{\Re}{\;\mathrm{Re}\;}
\renewcommand{\Im}{\;\mathrm{Im}\;}

\title[Asymptotic lower bounds]%
{{Asymptotic lower bounds for
a class of Schr\"odinger equations}}
\author{}

\author[Luis Vega and Nicola Visciglia]{{Luis Vega\\
Universidad del Pais Vasco, Apdo. 64\\
48080 Bilbao, Spain\\
email: mtpvegol@lg.ehu.es\\
tel.:++34-946015475, fax: ++34-946012516\\
\vspace{0.2cm}
and\\
Nicola Visciglia\\
Dipartimento di Matematica Universit\`a di Pisa\\
Largo B. Pontecorvo 5, 56100 Pisa, Italy\\
email: viscigli@dm.unipi.it\\
tel.: ++39-0502212294, fax: ++39-0502213224}}

\maketitle

\date{}

\begin{Abstract}
We shall study the following initial value problem:
\begin{equation}\label{abs}{\bf i}\partial_t u - \Delta u + V(x) u=0, \hbox{ } 
(t, x) \in {\mathbf R} \times {\mathbf R}^n, 
\end{equation}
$$u(0)=f,$$
where $V(x)$ is a real 
short--range potential, whose radial derivative satisfies some supplementary assumptions.
More precisely we shall present a family of identities satisfied by the solutions
to \eqref{abs} that generalizes the ones proved in 
\cite{PL}  and \cite{VV} in the free case.
As a by--product of these identities we deduce some uniqueness results
for solutions to        \eqref{abs},
and a lower bound for the so called local smoothing  which becomes an identity in a precise asymptotic sense.
\end{Abstract}

\section{Introduction}

We shall study the following initial value  problem:
\begin{equation}\label{cau}
{\bf i}\partial_t u - \Delta u + V(x) u=0, \hbox{ } 
(t, x) \in {\mathbf R} \times {\mathbf R}^n,
\end{equation}
$$u(0)=f$$
under suitable assumptions on $V(x)$.

Let us recall that if $V(x) \in L^\infty({\mathbf R}^n)$
then the operator 
$$L^2({\mathbf R}^n) \supset H^2({\mathbf R}^n)
\ni u\rightarrow -\Delta u + V(x) u \in L^2({\mathbf R}^n)$$
is self--adjoint (see \cite{RS} for the proof of this fact).
In particular one can apply
the classical Stone theorem in order to deduce the existence of a unique solution 
$u(t, x)\in {\mathcal C}_t(L^2({\mathbf R}^n))$ 
(here and in the sequel we shall denote by ${\mathcal C}_t(X)$
the space of continuous functions of one variable valued in the Banach space $X$)
to the Cauchy problem \eqref{cau}, provided that $f\in L^2({\mathbf R}^n)$.

Hereafter we shall denote by
$e^{{\bf i}t\Delta_V}f$
the unique solution to \eqref{cau} at time $t\in \mathbf R$.
Let us recall  that the following conservation law is satisfied:
\begin{equation}\label{charge}
\|e^{{\bf i}t\Delta_V}f\|_{L^2({\mathbf R}^n)}\equiv \|f\|_{L^2({\mathbf R}^n)}
\hbox{ } \forall t\in \mathbf R .
\end{equation}

Notice that this identity implies that the operators $e^{{\bf i}t\Delta_V}$
define a family of isometries on $L^2({\mathbf R}^n)$.
Moreover, as a by--product of the Stone theorem one can deduce the following implication:
\begin{equation}\label{stonereed}
f\in H^2({\mathbf R}^n) \Rightarrow e^{{\bf i}t \Delta} f \in {\mathcal C}^1_t (H^2({\mathbf R}^n))
\end{equation}
(here we have denoted by ${\mathcal C}_t^1(X)$
the space of functions of one variable valued in the Banach space $X$
with a continuous derivative).
It is also well--known that the following conservation law holds:
\begin{equation}\label{perturbedenergy}
\int_{{\mathbf R}^n} (|\nabla_x u(t, x)|^2 + V(x) |u(t, x)|^2) \hbox{ } dx
\end{equation}
$$=\int_{{\mathbf R}^n} (|\nabla_x f(x)|^2 + V(x) |f(x)|^2) \hbox{ }
dx \hbox{ } \forall t\in \mathbf R,
$$
and in particular 
\begin{equation}\label{energy}\int_{{\mathbf R}^n} |\nabla_x u(t, x)|^2 \hbox{ } dx\leq C
\int_{{\mathbf R}^n} (|\nabla_x f(x)|^2 + |f(x)|^2) \hbox{ }
dx \hbox{ } \forall t \in \mathbf R\end{equation}
provided that $V(x)\geq 0$ and $V(x)\in L^\infty({\mathbf R}^n)$.

In the sequel we shall assume that $V(x)$
satisfies the following decay assumption:
\begin{equation}\label{SR0}
0\leq V(x)\leq \frac C{(1+ |x|)^{1+\epsilon}}\hbox{ } 
\forall x\in {\mathbf R}^n
\end{equation}
where  
$\epsilon, C>0.$

We shall also assume either that $V(x)$
is decreasing in the radial variable, i.e.
\begin{equation}\label{decay}\partial_{|x|} V\leq 0,\end{equation}
or that
\begin{equation}\label{new}
\lim_{|x|\rightarrow \infty}|x|\partial_{|x|}V(x)=0.
\end{equation}

We shall specify in every theorem which kind of assumptions we assume
on the derivative of $V$.

In order to state our results let us introduce the perturbed Sobolev spaces
$\dot H^s_V({\mathbf R}^n)$,
whose norm is defined as follows:
\begin{equation}\label{sobolev}
\|f\|_{\dot H^s_V({\mathbf R}^n)}\equiv \left \|\left 
(\sqrt {-\Delta + V}\right )^s u\right \|_{L^2({\mathbf R}^n)}
\hbox{ } \forall s\geq 0.
\end{equation}

Our first result contains a family of identities satisfied by solutions to \eqref{cau}.
Let us underline that these identities 
represent a generalization of the identities proved in the free case, i.e. $V(x)\equiv 0$,
in \cite{PL} and \cite{VV}.
In fact in \cite{VV2} a similar family of identities has been proved for the solutions to the 
conformally invariant nonlinear Schr\"odinger equation.
Next we shall denote by $D^2 \psi$ the Hessian matrix of the function $\psi$.
\begin{thm}\label{mainUCSL0}
Let $u(t,x)$ be the solution to 
\eqref{cau} where $n\geq 1$, $f\in C^\infty_0({\mathbf R}^n)$.
Assume moreover that $V(x)$ 
satisfies \eqref{SR0} and one of the conditions 
\eqref{decay} or \eqref{new}.
Let $\psi$
be a radially symmetric function such that
the following limit exists
\begin{equation}
\label{limitder0}\lim_{|x|\rightarrow \infty} 
\partial_{|x|} \psi=\psi'(\infty)\in [0, \infty)
\end{equation}
and moreover
$$\nabla \psi, D^2 \psi, \Delta^2 \psi\in L^\infty({\mathbf R}^n)$$ 
Then the 
following identity holds:
\begin{equation}\label{vai}
\lim_{T\rightarrow \infty} \int_{-T}^T \int_{{\mathbf R}^n}  \left[\nabla_x \bar u D^2 \psi \nabla_x u
-
(\Delta^2 \psi +
4\partial_{|x|} V  \partial_{|x|}\psi )\frac{|u|^2}4
\right ]dx dt 
\end{equation}\begin{equation*}
= \psi'(\infty)\|f\|_{\dot H^\frac 12_V({\mathbf R}^n)}^2.
\end{equation*}
\end{thm}

As a by--product of the argument involved in the proof of theorem \ref{mainUCSL0}
we can construct a natural Banach space $\Sigma^\frac 12$ (whose definition will be 
given below)
that is invariant along
the flow associated to \eqref{cau}.
Moreover we shall deduce one uniqueness result for solutions to \eqref{cau}
provided that $f\in \Sigma^\frac 12$.
In order to define the space $\Sigma^\frac 12$
we first introduce the weighted Lebesgue space 
$L^2_{|x|}({\mathbf R}^n)$ defined as the completion of $C^\infty_0({\mathbf R}^n)$
with respect to the following norm:
\begin{equation}\label{weighted}
\|f\|_{L^2_{|x|}}^2 \equiv\int_{{\mathbf R}^n} |x| |f(x)|^2 \hbox{ } dx.
\end{equation}
The Banach space $\Sigma^{\frac 12}$
is defined as follows:
\begin{equation}\label{sigma}
\Sigma^\frac 12\equiv \dot H^\frac 12_V \cap L^2_{|x|}
\end{equation}
and can be endowed with the norm 
$$\|f\|_{\Sigma^\frac 12}^2\equiv  \|f\|_{\dot H^\frac 12_V}^2+
\|f\|_{L^2_{|x|}}^2.$$  

We can state our second result.

\begin{thm}\label{invariance}
Let $u(t,x)$ be the solution to 
\eqref{cau} where $n\geq 2$, $f\in \Sigma^\frac 12$ and
$V(x)$ 
satisfies the same assumptions as in theorem \ref{mainUCSL0}.
Then we have the following a--priori estimate:
\begin{equation}\label{famboc}\|u(t)\|_{\Sigma^\frac 12}^2\leq 
\|f\|_{L^2_{|x|}}^2 
+ C(1+ |t|) \|f\|^2_{\dot H^\frac 12_V}
\hbox{ } \forall t\in \mathbf R,
\end{equation}
for a suitable $C>0$.
\noindent In particular for every $t\in \mathbf R$
we have that $e^{{\bf i}t\Delta_V} f\in \Sigma^\frac 12$
provided that $f\in \Sigma^\frac 12$.
\noindent Moreover 
\begin{equation}\label{secunda}
\lim_{t\rightarrow \pm \infty} 
\int_{{\mathbf R}^n} \frac{|x|}{|t|} |u(t,x)|^2\hbox{ } dx
=2 \|f\|_{\dot H^\frac 12_V}^2.
\end{equation}
In particular if
$$\lim_{t\rightarrow \pm \infty} 
\int_{{\mathbf R}^n} \frac{|x|}{|t|} |u(t, x)|^2\hbox{ } dx=0$$
then $u\equiv 0$.
\end{thm}  

\begin{remark}
From a technical point of view we assume $n\geq 2$ in theorem \ref{invariance},
since the proof of lemma \ref{kpv} (that in turn is needed in the proof of 
theorem \ref{invariance})
does not work in dimension $n=1$.
\end{remark}

\begin{remark}
Along the proof of theorems \ref{mainUCSL0}
and \ref{invariance},
we shall make extensively use of 
the existence and completeness of the wave operator
under the assumptions \eqref{SR0} and \eqref{new} on $V(x)$
(see section \ref{agmsca}).
\end{remark}

\begin{remark}
In order to prove theorems 
\ref{mainUCSL0} and \ref{invariance}
we shall need some intemediate results,
whose proof in some cases
could be deduced by avoiding the use of the existence and completeness
of the wave operator. For instance
lemma \ref{anco} in section \ref{agmsca} 
follows from the general RAGE
theory (see \cite{RS2}). However we have proposed
a proof that involves the existence and completeness of the wave operator
in order to make the paper selfcontained as much as possible. 
In the appendix \ref{appendix}
we shall make some connections between
the classical RAGE theorem and our results.
\end{remark}

Next we shall deduce some direct consequences from the identity \eqref{vai}.
In particular we shall show how it allows us to prove 
a lower bound to the classical local smoothing estimate. For a proof of the local smoothing estimate
in the free case see \cite{CS}, \cite{Sj},
\cite{Vega} and also their extensions in \cite{bk}, \cite{ky},
\cite{S}. In particular in \cite{S}
the issue of the best constants involved 
in the local smoothing estimate is considered. 
 
First we shall present our results in dimension $n\geq 4$.

\begin{thm}\label{mainSL}
Let $u(t,x)$ be the solution to \eqref{cau},
where $n\geq 4$ and $V(x)$ satisfies \eqref{SR0}
and \eqref{decay}.
Then the following a--priori estimate is satisfied
for every $f\in \dot H^\frac 12_V({\mathbf R}^n)$:
\begin{equation}\label{eqSL}
 \|f\|_{\dot H^\frac 12_V({\mathbf R}^n)}^2 
\leq \sup_{R>0}\frac 1R \int_{-\infty}^\infty \int_{|x|<R}
|\nabla_x u|^2 dxdt \leq C \|f\|_{\dot H^\frac 12_V({\mathbf R}^n)}^2,\end{equation}
where $C>0$ is a suitable constant independent of $f$.
\end{thm}

In next result we give a better lower bound than the one in \eqref{eqSL}.

\begin{thm}\label{mainUCSL}
Let $u(t,x)$ be 
the solution to \eqref{cau} where $n\geq 4$, $f$ and 
$V(x)$ are as in theorem \ref{mainSL}.
Then we have:
\begin{equation}\label{eqUCSL}
\lim_{R\rightarrow \infty} \frac 1R 
\int_{-\infty}^\infty \int_{|x|<R}|\partial_{|x|} u|^2 dxdt
=\|f\|_{\dot H^\frac 12_V({\mathbf R}^n)}^2.
\end{equation}
Therefore if 
$$\lim_{R\rightarrow \infty} \frac 1R \int_{-\infty}^\infty 
\int_{|x|<R}|\partial_{|x|} u|^2 dxdt
=0$$
then $u\equiv0$.
\end{thm}

\begin{remark}\label{rembasic}
Notice that if you choose  formally $\psi\equiv |x|$ in \eqref{vai}
and if you work in dimension $n\geq 4$,
then the general identity \eqref{vai} becomes 
\begin{equation}\label{morawetz}
\int_{-\infty}^\infty \int_{{\mathbf R}^n} \left (\frac{|\nabla_\tau u|^2}{|x|}
+ \frac{(n-1)(n-3)}{4}\frac{|u|^2}{|x|^3} - 
\partial_{|x|}V |u|^2 \right )\hbox{ } dxdt
= \|f\|_{\dot H^\frac 12_V({\mathbf R}^n)}^2
\end{equation}
(here $\nabla_\tau  u$ denotes the angular part
of the gradient of $u$). 
Let us underline that in the case $V(x)\equiv 0$
the previous identity has been proved in \cite{PL} with a different approach.
Moreover \eqref{morawetz} represents a precised version
of the result in \cite{lis} and \cite{m}, where 
\eqref{morawetz} is stated as an inequality 
and not as an identity.
Notice also that the function $\psi \equiv |x|$ does not satisfies all the assumptions
required in theorem \ref{mainUCSL0}.
However in order to make precise the argument involved 
in the proof of \eqref{morawetz}, it is sufficient to choose in \eqref{vai}
the test function $\psi$ to be equal to
$\sqrt{\epsilon^2 + |x|^2}$ and to get the limit 
in the corresponding identity as $\epsilon \rightarrow 0$.
An alternative way to prove properly \eqref{morawetz} 
it is to combine the proof of \eqref{vai}
with the argument used in \cite{m} (in fact in \cite{m}
the integration by parts technique that
we use in the proof of \eqref{vai} is completely 
justified also when $\psi\equiv |x|$).
\end{remark}

\begin{remark}
Let us point out that theorems
\ref{mainSL} and \ref{mainUCSL} are stated in dimension $n\geq 4$, while theorem \ref{mainUCSL0} is stated in any dimension $n\geq 1$. 
The main reason is that
in order to take advantage of the identity \eqref{vai}
we shall choose the test function $\psi(x)$
in a suitable way. 
As it will be clear in the sequel,
we shall be able to make such a good choice
only in dimension $n\geq 4$.
\end{remark}

In dimension $n=3$ we are able
to prove the following result.

\begin{thm}\label{mainSL3}
Let $u(t,x)$ be the solution to \eqref{cau}
where $n=3$, $V(x)$ satisfies \eqref{SR0} 
and \eqref{decay}.
Then the following a priori estimate is satisfied
for every $f\in \dot H^\frac 12_V({\mathbf R}^3)$:
\begin{equation}\label{eqSL3}
 c \|f\|_{\dot H^\frac 12_V({\mathbf R}^3)}^2 
\leq \sup_{R>0}\frac 1R \int_{-\infty}^\infty \int_{|x|<R}
\left ( |\nabla_x u|^2  
+ \frac 1{R^2} 
|u|^2 \right ) dxdt 
\leq C \|f\|_{\dot H^\frac 12_V({\mathbf R}^3)}^2,\end{equation}
where $c, C>0$ are constants independent of $f$.
\end{thm}

\begin{thm}\label{mainUCSL3}
Let $u(t,x)$ be the
solution to \eqref{cau} where $n=3$, $f$ and 
$V(x)$ are as in theorem \ref{mainSL3}.
Then we have:
\begin{equation}\label{eqUCSL3}
\liminf_{R\rightarrow \infty} \frac 1R 
\int_{-\infty}^\infty \int_{|x|<R} \left( |\partial_{|x|} u|^2 dxdt
+ \frac 1{R^2}
|u|^2 \right ) dxdt
\geq \frac 12  \|f\|_{\dot H^\frac 12_V({\mathbf R}^3)}^2.
\end{equation}
Therefore if 
$$\liminf_{R\rightarrow \infty} \frac 1R \int_{-\infty}^\infty 
\int_{|x|<R}\left (|\partial_{|x|} u|^2 
+ \frac 1{R^2} 
|u|^2 \right ) dxdt=0$$
then $u\equiv0$.
\end{thm}

\begin{remark}
Starting with  \eqref{vai} it is possible to show
the following version of the identity \eqref{morawetz}
in dimension $n=3$:

\begin{equation}\label{morawetz3}
\int_{-\infty}^\infty \int_{{\mathbf R}^3} \left (\frac{|\nabla_\tau u|^2}{|x|}
-\partial_{|x|}V |u|^2 \right )\hbox{ } dxdt
+ \frac 32 \pi \int_{-\infty}^\infty
|u(0, t)|^2 dt = \|f\|_{\dot H^\frac 12_V({\mathbf R}^3)}^2.
\end{equation} 

Exactly as for \eqref{morawetz} the previous identity 
has been proved previously in \cite{PL} in the free case 
and it represents a precised version
of a result proved in \cite{lis} and \cite{m}, where 
\eqref{morawetz3} is stated as an inequality.

Indeed the proof of \eqref{morawetz3} follows formally by choosing the function $\psi \equiv |x|$ 
in \eqref{vai}. However in dimension $n=3$ the function $\psi\equiv |x|$ is very singular
since its bilaplacian is a multiple of the Dirac delta
and hence in order to justify all the computations we have to argue
as in dimension $n\geq 4$ (see remark \ref{rembasic}). 
\end{remark}

The paper is organized as follows.
In section \ref{agmsca} we shall prove 
some asymptotic properties of solutions to \eqref{cau}.
Sections \ref{BAR} and \ref{INV}
will be devoted to the proof of theorems
\ref{mainUCSL0} and \ref{invariance}.
The proof of theorems \ref{mainSL} and \ref{mainUCSL}
will be given in section
\ref{benfra}, while theorems \ref{mainSL3} and \ref{mainUCSL3}
will be proved in section \ref{fanucci}.
Finally in the appendix \ref{appendix} 
we shall discuss
some connections between our results and the classical RAGE theorem.
\vspace{0.1cm}

Next we shall fix some notations.

\vspace{0.1cm}

\noindent {\bf Notations.}
For every potential $V(x) \geq 0$ and for every real number $s\geq 0$ we shall denote by 
$\dot H^s_V$ the perturbed Sobolev space whose norm is defined
in \eqref{sobolev}.
In particular when $V\equiv 0$ these spaces reduce to the standard Sobolev spaces
$\dot H^s$ whose norm is defined as follows
$$\|f\|_{\dot H^s}^2\equiv \int_{{\mathbf R}^n} |\hat f(\xi )|^2|\xi|^{2s}
d\xi,$$
where 
$$\hat f(\xi)\equiv \int_{{\mathbf R}^n}e^{-2\pi 
{\bf i} x\xi} f(x) dx.$$
\noindent In some cases we shall also write
$${\mathcal F }(f)\equiv \hat f.$$

The spaces $L^2_{|x|}$ and $\Sigma^\frac 12$
are respectively the ones defined in \eqref{weighted} and \eqref{sigma}.

For any $1\leq p, q\leq \infty$ 
$$L^p_x \hbox{ and } L^p_t L^q_x$$
denote the Banach spaces
$$L^p({\mathbf R}^n) \hbox{ and } L^p({\mathbf R}; L^q({\mathbf R}^n)).$$
We shall also write
$$L^p_t L^p_x\equiv  L^p_{t, x}.$$

For every $V(x)\in L^\infty_x$ 
we shall denote by $e^{{\bf i }t\Delta_V}$
the group associated to \eqref{cau} via the Stone theorem.

Given any couple of Banach spaces $X$ and $Y$,
we shall denote by ${\mathcal L}(X, Y)$
the space of linear and continuous
functionals between $X$ and $Y$.
 
Given a space--time dependent function $w(t, x)$ we shall denote by
$w(t_0)$ the trace of $w$ at fixed 
time $t\equiv t_0$, in case that it is 
well--defined.

We shall denote by $\int ...\hbox {  }dx ,
\int ... \hbox{ } dt$
and $\int \int ...\hbox{ } dx dt$ the 
integral of suitable functions
with respect to  space, time, 
and space--time variables respectively.

When it is not better specified we shall 
denote by $\nabla v$ the gradient
of any time--dependent function $v(t, x)$ 
with respect to the space variables.
Moreover $\nabla_\tau$ and $\partial_{|x|}$
shall denote respectively the angular 
gradient and the radial derivative.

If $\psi\in C^2({\mathbf R}^n)$, then 
$D^2 \psi$ will represent the hessian matrix of $\psi$.

Given a set $A\subset {\mathbf R}^n$ we denote by $\chi_A$
its characteristic function. 

We shall use the function
$$\langle x \rangle\equiv \sqrt{1+ |x|^2}.$$

\section{Wave operators and asymptotic behaviour of solutions}
\label{agmsca}

Let us recall that if $V(x)$ satisfies  \eqref{SR0},
then the wave operators ${\mathcal W}_\pm$ are well--defined and complete
(see \cite{A}, \cite{Enss}, \cite{RS2} and \cite{RS4}). More precisely
for every $f\in L^2_x$ there exist two functions ${\mathcal W}_\pm(f)\in L^2_x$
uniquely defined and such that
\begin{equation}\label{scattering}
\lim_{t\rightarrow \pm \infty} \|u(t) - e^{{\bf i}t\Delta} {\mathcal W}_\pm (f)\|_{L^2_x}
=0,
\end{equation}
where $u\in {\mathcal C}_t(L^2_x)$ denotes the unique solution to
\eqref{cau} with initial data $f$ and $e^{{\bf i}t \Delta}$
represents the propagator at time $t$ associated to the 
free Schr\"odinger equation, i.e. \eqref{cau}
with $V(x)\equiv 0$ (for a proof of \eqref{scattering} see \cite{A}).

\vspace{0.1cm}

In the sequel we shall need 
the following asymptotic description of the free waves,
whose proof can be found in \cite{RS}. 
\begin{prop}
Assume $f\in L^2_x$ and $n\geq 1$, then:
\begin{equation}\label{asymptotique+}
\lim_{t\rightarrow \pm \infty} \left \|e^{{\bf i}t\Delta }f - 
e^{\mp{\bf i}n\pi/4}\frac{e^{\pm {\bf i} \frac{|x|^2}
{4t}}}{(4\pi t)^{n/2}} 
\hat f \left(\pm \frac { x}{4\pi t} \right)
\right\|_{L^2_x}=0.
\end{equation}
\end{prop}

Since now on we shall denote by ${\mathcal W}_\pm$ the wave operators 
defined above and by $\mathcal F$
the Fourier transform.

Next we shall state one of the basic results of this
paper.

\begin{prop}\label{scas2}
Let $u(t, x)$ be the solution to \eqref{cau} where $n\geq 1$, $f\in C^\infty_0({\mathbf R}^n)$
and $V(x)$ satisfies \eqref{SR0} 
and one of the two conditions \eqref{decay} or \eqref{new}.
Assume that
$\psi$
is a radially symmetric such that the following limit exists:
\begin{equation}
\label{limitder}\lim_{|x|\rightarrow \infty} \partial_{|x|} \psi=\psi'(\infty)\in [0, \infty).
\end{equation}
\noindent Then
\begin{equation}\label{somcon77}\lim_{t\rightarrow \pm \infty}  {\mathcal Im}
\left (\int \bar u(t) \nabla u(t) \cdot \nabla \psi \hbox{ } dx \right ) 
= \mp 2\pi \psi'(\infty)\int |x| |g_\pm(x)|^2\hbox{ } dx\end{equation}
where
$g_\pm={\mathcal F}[{\mathcal W}_\pm f]$.
Moreover the following identity holds:
\begin{equation}\label{dot7}
\|f\|_{\dot H^\frac 12_V}^2= 2\pi \int |x||g_\pm(x)|^2 \hbox{ } dx.
\end{equation}
\end{prop}

We shall need the following lemma which is a consequence of the RAGE theorem. 
In order to be self--contained,
we have decided to include
a proof of it which is
based on the existence and completeness of the wave operators
${\mathcal W}_\pm$ introduced at the beginning of the section.

\begin{lem}\label{anco}
Let $u(t, x)$ be the solution to \eqref{cau} where
$n\geq 1$, $f\in L^2_x$ and $V(x)$ satisfies \eqref{SR0},
then
\begin{equation}\label{anco3}
\lim_{t\rightarrow \pm \infty} \frac 1t \int_0^t \int W(x)|u|^2 \hbox{ } dxds =0
\end{equation}
where $W\in L^\infty_x$ is such that
$$\lim_{|x|\rightarrow \infty} W(x)=0.$$
\end{lem}

\noindent{\bf Proof.} 
We shall prove the result only as $t\rightarrow \infty$, since
the case $t\rightarrow -\infty$ can be treated similarly.

Notice the following identity
\begin{equation}\label{anco5}\int_0^t \int W(x)|u|^2 \hbox{ } dxds
=I(t)+ II(t) 
\hbox{ } \forall t>0
\end{equation}
where
$$I(t)=\int_0^t\int W(x)\left[|u|^2 - 
\frac{1}{(4\pi s)^n} \left |\hat g \left(\frac x{4\pi s} \right)\right |^2\right ] \hbox{ } dxds
$$
and
$$
II(t)=\int_0^t\int W(x)\left |\hat 
g\left (\frac x{4 \pi s}\right )\right |^2 \hbox{ } \frac{dxds}{(4\pi s)^n}.
$$
where $g\equiv {\mathcal W}_+(f)$. Due to \eqref{asymptotique+}
and to the assumption $W\in L^\infty_x$ we can deduce 
\begin{equation}\label{rauch1}
I(t) \leq \|W\|_{L^\infty_x}\int_0^t h(s) ds \hbox{ and } \lim_{s\rightarrow \infty} h(s)=0
\end{equation}
where $$h(s)= \int \left[|u(s)|^2 - 
\frac{1}{(4\pi s)^n} \left |\hat g \left(\frac x{4\pi s} \right)\right |^2\right ] \hbox{ } dx
.$$

On the other hand we have
$$II(t)=\int_0^t\int W(4\pi sx) |\hat g(x)|^2 \hbox{ } ds dx 
$$
that due to the dominated convergence theorem and to the decay assumption
made on $W(x)$ implies
\begin{equation}\label{rauch2}
II(t)=\int_0^t H(s) ds \hbox{ and } \lim_{s\rightarrow 0} H(s)=0
\end{equation}
where
$$H(s)=\int W(4\pi sx) |\hat g(x)|^2 \hbox{ } dx.$$

By combining \eqref{anco5} with \eqref{rauch1} and \eqref{rauch2}
it is easy to deduce \eqref{anco3}.

\hfill$\Box$

\begin{remark}
Notice that in order to deduce
\eqref{anco3} we have shown
that
\begin{equation}\label{RAGE}
\lim_{t\rightarrow \pm \infty}
\int W(x)|u(t)|^2 \hbox{ } dx=0
\end{equation}
which is stronger than  \eqref{anco3}.
In fact \eqref{anco3} could be proved for a much larger class
of potentials $V(x)$ by using the general RAGE theorem.
In the appendix \ref{appendix}
we shall show how 
to deduce \eqref{RAGE} by using as a starting point 
\eqref{anco3}. 
\end{remark}

\begin{lem}\label{ragecazenave}
Let $u(t, x)$ be the solution to \eqref{cau} where $n\geq 1$, $f\in C^\infty_0({\mathbf R}^n)$ 
and $V(x)$ satisfies \eqref{SR0} and 
one of the two conditions \eqref{decay}
or \eqref{new}. Then
we have:
\begin{equation}\label{or}
\lim_{t\rightarrow \pm \infty} \left \|\frac xt u(t) - 2{\bf i} \nabla u(t)
\right \|_{L^2_x}=0.
\end{equation} 
\end{lem}

\noindent {\bf Proof.} We prove \eqref{or} only for
$t\rightarrow \infty$, since the case $t\rightarrow -\infty$ is similar.\\

\vspace{0.1cm}

{\em First case: $V(x)$ satisfies \eqref{decay}}

\vspace{0.2cm}

\noindent The following identity is well--known (see \cite{caze}):
\begin{equation}\label{conservation}
\|x u(t) - 2{\bf i}t \nabla u(t)\|_{L^2_x}^2 + 
 4 t^2 \int V(x) |u(t)|^2 dx 
\end{equation}
\begin{equation*}
= \int |x|^2 |f(x)|^2 \hbox{ } dx 
+ \int_0^t s \theta(s) \hbox{ } ds \hbox{  } \forall t\in \mathbf R
\end{equation*}
where
$$\theta(s)= 8 \int \left
(V(x) + \frac 12 |x| \partial_{|x|} V(x)\right ) |u(s)|^2 \hbox{ } dx.$$

By combining the sign assumption done on $V(x)$ and $\partial_{|x|} V$ 
with \eqref{conservation} we get:
\begin{equation}\label{ragecaze1}
\left \|\frac xt u(t) - 2{\bf i} \nabla u(t)
\right \|_{L^2_x}^2 \end{equation}
$$\leq \frac{\int |x|^2 |f(x)|^2  \hbox{ } dx}{t^2}
+\frac{8 \int_0^t \left(\int V(x) |u(s)|^2\hbox{ } dx \right ) \hbox{ } ds}{t}.$$

By combining this inequality with \eqref{anco3} we get \eqref{or}.

\vspace{0.1cm}

{\em Second case: $V(x)$ satisfies \eqref{new}}

\vspace{0.2cm}

\noindent In this case
we can use \eqref{conservation} as above and we can deduce 
\begin{equation}\label{ragecaze}
\left \|\frac xt u(t) - 2{\bf i} \nabla u(t)
\right \|_{L^2_x}^2 \leq \frac{\int |x|^2 |f(x)|^2 \hbox{ } dx}{t^2}
+\frac{\int_0^t \left(\int W(x) |u(s)|^2\hbox{ } dx \right ) \hbox{ } ds}{t}
\end{equation}
where 
$$W(x)=8\left(V(x) + \frac 12 |x| \partial_{|x|}V(x)\right ).$$

Notice that due to \eqref{SR0} and \eqref{new} we have that $\lim_{|x|\rightarrow \infty} W(x)=0$,
then we can  use \eqref{anco3} in order to deduce \eqref{or}.

\hfill$\Box$

\begin{remark}\label{Enssremark}
Let us underline that in \cite{Enss}
it is proved the 
existence of a sequence $\{t_n\}_{n\in \mathbf N}$
such that:
\begin{equation}\label{enssweak}\lim_{n\rightarrow \infty} t_n=\infty
\hbox{ and } \lim_{n\rightarrow \infty} \left \|\frac x{t_n} u(t_n) - 
2{\bf i} \nabla u(t_n)
\right \|_{L^2_x}=0.\end{equation}
Notice that \eqref{enssweak}
is a weaker version of \eqref{or},
however in \cite{Enss} it is a basic tool
in order to prove the completeness of the wave operators,
provided that $V(x)$ satisfies 
the assumptions \eqref{SR0} and \eqref{new}.
\end{remark}
\begin{remark}
Notice that we obtain \eqref{or}
using lemma \ref{anco} which we prove using
the completeness of the wave operators.
In the appendix \ref{appendix}
we shall give a proof of
\eqref{or} that does not involve
a--priori the completeness of the wave operator.
Moreover we shall show that 
\eqref{or} is still satisfied 
for a class of potentials
$V(x)$ more general than the ones
that satisfy the decay
assumptions \eqref{SR0} and \eqref{new}.
\end{remark}

\noindent {\bf Proof of proposition \ref{scas2}}
As usual we treat only the case $t\rightarrow \infty$,
the case $t\rightarrow -\infty$ can be treated similarly.

Along the proof we shall use the function $g(x)$ defined as
follows
$$g\equiv {\mathcal F} [{\mathcal W}_+ (f)].$$

Let us introduce the following identity:
\begin{equation}\label{splittingR}
{\mathcal Im}\int \bar u(t) \nabla u(t) \cdot \nabla \psi  \hbox{ } dx
=I(t, R)+ II(t, R) \hbox{ } \forall t\in {\mathbf R}, R>0,
\end{equation}
where
$$I(t, R)={\mathcal Im}\left (
\int_{|x|>4\pi Rt} \bar u(t) \nabla u(t) \cdot \nabla \psi  \hbox{  } dx\right )
$$ and 
$$II(t, R)={\mathcal Im}\left (\int_{|x|<4\pi Rt} \bar u(t) \nabla u(t) \cdot \nabla \psi \hbox{ } dx\right ).$$
 
\vspace{0.1cm}

{\em Estimate for $I(t, R)$}

\vspace{0.2cm}

\noindent Notice that the Cauchy--Schwartz inequality
implies:
\begin{equation}\label{cs}\left |\int_{|x|>4\pi R t} \bar u(t) \nabla u(t) 
\cdot \nabla \psi \hbox{ } dx \right |
\leq  C \|\nabla \psi\|_{L^\infty_x} \|f\|_{H^1_x} 
\left (\int_{|x|>4 \pi Rt} |u(t)|^2 dx\right )^\frac 12.
\end{equation}
where we have used \eqref{energy}.

On the other hand due to \eqref{asymptotique+} and due to the definition of $g(x)$ we get:
$$\lim_{t\rightarrow \infty}
\int_{|x|>4\pi Rt} \left [|u(t)|^2 - \frac{1}{(4\pi t)^{n}} \left|
g\left( \frac x {4\pi t}\right)\right|^2\right ]dx
=0$$
and then
\begin{equation}\label{lim}\lim_{t\rightarrow \infty}
\int_{|x|>4\pi Rt} |u(t)|^2 dx = \int_{|x|>R} |g(x)|^2dx.\end{equation}

Since $g\in L^2_x$ we can combine \eqref{cs}
with \eqref{lim} in order to deduce that
\begin{equation}\label{I}
\forall \epsilon>0 \hbox{ } \exists R(\epsilon)>0
\hbox{ s. t.}
\limsup_{t\rightarrow \infty} |I(t, R)|<\epsilon \hbox{ } \forall R>R(\epsilon).
\end{equation}

\vspace{0.1cm}

{\em Estimate for $II(t, R)$}

\vspace{0.2cm}

\noindent Notice that \eqref{charge} and \eqref{or} imply:
\begin{equation}\label{som2}\lim_{t\rightarrow \infty}
\left [\int_{|x|<4\pi Rt} \bar u(t) \nabla u(t) \cdot \nabla \psi \hbox{ } dx 
+\frac {\bf i} {2t} \int_{|x|<4\pi Rt} |x|   
\partial_{|x|}\psi  |u(t)|^2 \hbox{ } dx\right]=0.
\end{equation}

Moreover we have the following identities:
\begin{equation}\label{limintSL}  
\int_{|x|<4\pi Rt} |x| 
\partial_{|x|}\psi  |u(t)|^2 \hbox{ } \frac{dx}{t}
\end{equation}
$$=  \int_{|x|<4\pi Rt} |x| \partial_{|x|}\psi \left [|u(t)|^2 -\frac 1{(4\pi t)^n} 
\left | g\left (\frac x{4\pi t} \right )\right |^2\right ]
\hbox{ }\frac{dx}{t}
$$$$+ \frac 1{(4\pi)^n}\int_{|x|<4\pi Rt} |x| \left ((\partial_{|x|} \psi 
-\psi'(\infty)\right)  \left | g
\left ( \frac x {4\pi t}\right)\right |^2 \frac{dx}{t^{n+1}}
$$
$$
+ \frac 1{(4\pi)^n}\psi'(\infty) \int_{|x|<4\pi Rt}
|x| \left |g\left (\frac x{4\pi t} \right )
\right |^2 \frac{dx}{t^{n+1}}.
$$

Notice that the following estimate is trivial:
\begin{equation}\label{lim4SL}
\left | \int_{|x|<4\pi Rt}
|x| \partial_{|x|}\psi  \left [|u(t)|^2 -\frac 1{(4\pi t)^n} 
\left | g\left (\frac x{4\pi t} \right )\right |^2\right ]
\hbox{ } \frac{dx}{t}
\right | \end{equation}
\begin{equation*}
\leq 4\pi R \|\partial_{|x|}\psi\|_{L^\infty_x}\int
\left ||u(t)|^2 -\frac 1{(4\pi t)^n} 
\left | g\left (\frac x{4\pi t} \right )\right |^2\right |dx
\rightarrow 0 \hbox{ as } t\rightarrow \infty,
\end{equation*}
where at the last step we have combined \eqref{scattering}
with \eqref{asymptotique+}.

Moreover the change of variable formula implies:
$$\frac 1{(4\pi)^n}\left |\int_{|x|< 4\pi Rt} |x| \left ((\partial_{|x|} \psi 
-\psi'(\infty)\right ) \left | g\left ( \frac x 
{4 \pi t}\right)\right |^2 \frac{dx}{t^{n+1}}\right |
$$
$$\leq 4\pi R  \int \left | \partial_{|x|} \psi ( 4\pi t x )
-\psi'(\infty)\right| |g(x)|^2 dx,$$ 
that in conjunction  with 
the dominated convergence theorem and with assumption
\eqref{limitder} implies:
\begin{equation}
\label{lim7SL}\lim_{t\rightarrow \infty} \frac 1{(4\pi)^n}\int_{|x|< 4\pi Rt} 
|x| \left ((\partial_{|x|} \psi 
-\psi'(\infty) \right) 
\left | g\left ( \frac x {4\pi t}\right)\right |^2 \frac{dx}{t^{n+1}}
=0.\end{equation}

Due again to the change of variable formula we get
$$
\frac {\psi'(\infty) }{(4\pi)^n} \int_{|x|< 4\pi Rt} |x| 
\left |g\left (\frac x{4 \pi t} \right )\right |^2 \frac{dx}{t^{n+1}}
= 4\pi \psi'(\infty)\int_{|x|< R } |x| |g(x)|^2 dx,
$$
and in particular
\begin{equation}\label{lim14SL}
\lim_{t\rightarrow \infty}
\frac{\psi'(\infty)}{(4\pi)^n} \int_{|x|< 4\pi Rt} 
|x| \left |g\left (\frac x{4\pi t} \right )\right |^2 \frac{dx}{t^{n+1}}
= 4 \pi \psi'(\infty)\int_{|x|< R} |x| |g(x)|^2 \hbox{ } dx.
\end{equation}

By combining 
\eqref{lim4SL},\eqref{lim7SL}, \eqref{lim14SL} with \eqref{som2}
and \eqref{limintSL} we deduce
\begin{equation}\label{II}
\lim_{t\rightarrow \infty} II(t, R)=
-2\pi {\psi'(\infty)} \int_{|x|< R} |x| |g(x)|^2 dx.
\end{equation}

By combining \eqref{splittingR} with \eqref{I} and \eqref{II} we get
\eqref{somcon77} at least in the case $t\rightarrow \infty$.
The other case is similar.

\vspace{0.1cm}

{\em Proof of \eqref{dot7}}

\vspace{0.2cm}

\noindent Recall that ${\mathcal W}_+:L^2_x\rightarrow L^2_x$
is an isometry and moreover
$${\mathcal W}_+ \circ f(-\Delta)= f(-\Delta +V) \circ {\mathcal W}_+.$$

By combining these facts with the definition of $g$, i.e. $g\equiv {\mathcal F}
[{\mathcal W}_+ f]$, we get: 
$$\int |x||g(x)|^2 \hbox{ } dx=\|{\mathcal W}_+ f\|_{\dot H^\frac 12_x}^2=
\frac 1{2\pi}\|(-\Delta)^\frac 14 \circ {\mathcal W}_+f\|_{L^2_x}^2$$$$
=\frac 1{2\pi}\|{\mathcal W}_+ \circ (-\Delta_V)^\frac 14 f\|_{L^2_x}^2
=\frac 1{2\pi}\|(-\Delta_V)^\frac 14 f\|_{L^2_x}^2=
\frac 1{2\pi} \|f\|_{\dot H^\frac 12_{V}}^2.$$
 
\hfill$\Box$

\section{Proof of theorem
\ref{mainUCSL0} and some consequences}\label{BAR}

\noindent {\bf Proof of theorem \ref{mainUCSL0}}
Following
\cite{BRV} we  multiply
\eqref{cau} by the quantity 
\begin{equation}\label{multiplier}
\nabla \bar u  \cdot \nabla \psi+ \frac 12 \bar u \hbox{ } \Delta \psi ,
\end{equation}
and we integrate on the strip $(-T, T)\times \mathbf R^n$.
In this way we get the
following family of identities:
\begin{equation}\label{BarcRuiVeg}
\int_{-T}^T\int  \left(\nabla \bar u D^2 \psi \nabla u
-\Delta^2 \psi \frac{|u|^2}4-
\partial_{|x|} V  \partial_{|x|}\psi |u|^2 
\right )dx dt\end{equation}
\begin{equation*} = -\frac 12
{\mathcal Im}\sum_\pm
\int \bar u(\pm T) \nabla u(\pm T) \cdot \nabla \psi \hbox{ } dx,
\end{equation*}
(for more details on this computation see \cite{BRV} and \cite{VV}).

Indeed all the integration by parts 
involved in the proof of \eqref{BarcRuiVeg}
can be completely justified by a density argument due to \eqref{stonereed}.

Notice that the identity \eqref{vai} follows by combining 
\eqref{somcon77}, \eqref{dot7}
and \eqref{BarcRuiVeg}.

\hfill$\Box$

Next we shall exploit \eqref{vai} in order to deduce 
some  a--priori estimates satisfied by the solutions to \eqref{cau}.

\begin{lem}\label{brand}
Assume that $u(t, x)$ solves \eqref{cau} where $n\geq 4$,
$f\in \dot H^\frac 12_V$,
$V(x)$ satisfies \eqref{SR0}  and \eqref{decay},
then
\begin{equation}\label{bilaplacian}
\int\int \frac{1}{\langle x
\rangle^3} |u|^2 \hbox{ } dxdt<\infty
\end{equation}
and
\begin{equation}\label{radial}
\int\int |\partial_{|x|}V| |u|^2 \hbox{ }
dxdt<\infty. \end{equation}

\end{lem}

\noindent{\bf Proof.} Choose in \eqref{vai} the function $\psi(x)\equiv
|x|$ that is clearly a radially symmetric 
and convex function. Moreover we have 
$\partial_{|x|} \psi\equiv 1$ and
\begin{equation}\label{ngeq4}\Delta^2 (|x|) 
= -\frac{(n-1)(n-3)}{|x|^3}
 \hbox{ } \forall x\in {\mathbf R}^n
\hbox{ where }n\geq 4.
\end{equation}

Hence by choosing $\psi(x)\equiv |x|$
in \eqref{vai}, it is easy to deduce that
\begin{equation}\label{Bra}\frac{(n-1)(n-3)}4 
\int \int \frac{|u|^2}{|x|^3}\hbox{ } dxdt
- \int \int 
\partial_{|x|}V  |u|^2 \hbox{ }dx dt
\leq \|f\|_{\dot H^\frac 12_V}^2\end{equation}
where $C>0$ is a suitable constant,
and hence we get easily 
\eqref{bilaplacian} and \eqref{radial}.

\hfill$\Box$

Next we state a version of \eqref{radial} in dimension $n=3$.

\begin{lem}
Assume that $u(t, x)$ 
is solution to \eqref{cau} with $n=3$, $f\in \dot H^\frac 12_V$, 
$V(x)$ satisfies \eqref{SR0} 
and \eqref{decay}, 
then:
\begin{equation}\label{radial3}
\int\int |\partial_{|x|}V| |u|^2 \hbox{ }
dxdt<\infty.\end{equation}
\end{lem}

\noindent{\bf Proof.}
The proof of \eqref{radial3} is identical to the proof
of \eqref{radial}. Notice that 
by choosing in \eqref{vai} the test function $\psi(x)\equiv |x|$
and arguing as in the proof of
lemma \ref{brand} we
get
$$\frac 32 \pi \int
|u(0, t)|^2 dt -\partial_{|x|}V |u|^2 \hbox{ } dxdt
\leq \|f\|_{\dot H^\frac 12_V({\mathbf R}^3)}^2$$
where we have used
the property
\begin{equation}\label{ngeq3}-\Delta^2 (|x|)=6\pi \delta_0 \hbox{ on } {\mathbf R}^3.
\end{equation}

\hfill$\Box$

We can now deduce the following
\begin{prop}\label{noma}
Assume that $u(t, x)$ 
is a solution to \eqref{cau} with $n\geq 4$, $f\in \dot H^\frac 12_V$, 
$V(x)$ satisfies \eqref{SR0} 
and \eqref{decay}, then:
\begin{equation}\label{uno}
\lim_{R\rightarrow \infty}  \int \int 
 |\Delta^2 \phi_R| |u|^2
\hbox{ } dxdt =0,\end{equation}
where $\phi$ is a radially symmetric function
such that 
\begin{equation}\label{cil}|\Delta^2 \phi |\leq \frac{C}
{\langle x\rangle^3} \hbox{  } \forall x\in {\mathbf R}^n\end{equation}
and $\phi_R =R \phi\left (\frac{x}R\right )$.
\end{prop}

\noindent{\bf Proof.}
Notice that \eqref{cil} trivially implies
$$\int \int  |\Delta^2 \phi_R| |u|^2 \hbox{ } dx dt
\leq C  
\int\int \frac {|u|^2}{R^3 + |x|^3} \hbox{  } dxdt
\rightarrow 0 \hbox{ as } R\rightarrow \infty,$$
where we have combined 
the dominated convergence theorem with \eqref{bilaplacian}.

\hfill$\Box$

\begin{prop}\label{noma3}
Assume that $u(t, x)$ 
is solution to \eqref{cau} with $n\geq 3$, $f\in \dot H^\frac 12_V$, 
$V(x)$ satisfies \eqref{SR0} 
and \eqref{decay}, 
then:
\begin{equation}\label{due3}
\lim_{R\rightarrow \infty}  \int \int  
 |\partial_{|x|}V|
|\partial_{|x|} \phi_R | |u|^2 \hbox{ } dxdt=0
\end{equation}
where $\phi$ is a radially symmetric function
such that 
$$\partial_{|x|} \phi(0)=0,
|\partial_{|x|}\phi| \leq C \hbox{  } \forall x\in {\mathbf R}^n$$
and $\phi_R =R \phi\left (\frac{x}R\right )$.
\end{prop}

\noindent {\bf Proof.} 
It follows from the following identity:
\begin{equation}\label{largevar} \int\int 
 |\partial_{|x|}V|
|\partial_{|x|} \phi_R ||u|^2 \hbox{ } dxdt
\end{equation}
$$
=\int
\int |\partial_{|x|}V| \left|\partial_{|x|}\phi\left ( \frac xR\right)
\right ||u|^2 \hbox{ } dxdt
\rightarrow 0 \hbox{ as }
R\rightarrow \infty,$$
where at the last step we have combined the dominated convergence theorem
with \eqref{radial} (or with \eqref{radial3} in the specific case $n=3$) 
and with the assumption $\partial_{|x|}\phi(0)=0$.

\hfill$\Box$ 
 
\section{Proof of theorem \ref{invariance}}
\label{INV}
We shall need the following lemma, whose proof in dimension $n\geq 3$
follows an argument in \cite{BRV}.
\begin{lem}\label{kpv}
Assume that $n\geq 2$ and
$h\in L^2_x \cap \dot H^1_x$,
then we have the following inequality:
\begin{equation}\label{furt}
\int \bar h(x) \hbox{ } \nabla 
h(x) \cdot \frac{x}{|x|} \hbox{ } dx \leq C\|h\|_{\dot H^\frac 12_x}^2,
\end{equation}
where $C>0$ is a constant that depends on $n$.
\end{lem}

\noindent{\bf Proof.}
We introduce the quadratic form
$$a(f, g)=  \int \bar f(x) \hbox{ } \nabla g(x) 
\cdot \frac{x}{|x|} \hbox{ } dx.$$
Notice that the Cauchy--Schwartz inequality implies:
\begin{equation}\label{bil1}
|a(f, g)|\leq \|f\|_{L^2_x} \|g\|_{\dot H^1_x}.
\end{equation}
Next we split the proof in two cases.

\vspace{0.1cm}

{\em First case: $n\geq 3$}

\vspace{0.2cm}
\noindent 
By assuming $f, g$ regular enough we can use integration by parts
in order to deduce:
$$a(f, g)= - \int g(x) \nabla \bar f(x) \cdot \frac x{|x|} \hbox{ } dx
- (n-1)\int \frac 1{|x|}\bar f (x) g(x) \hbox{ } dx.$$

By combining the Hardy inequality and the Cauchy--Schwartz inequality with
the previous identity we get:
\begin{equation}\label{bil2}
|a(f, g)|\leq C \|f\|_{\dot H^1_x} \|g\|_{L^2_x}
\end{equation}

Notice that \eqref{furt} will follow by interpolation from
\eqref{bil1} and \eqref{bil2} and by choosing $f=g=h$.

\vspace{0.1cm}

{\em Second case: $n=2$}

\vspace{0.2cm}

\noindent 
In this case we are not allowed to use
the Hardy inequality in order to deduce \eqref{bil2}.
Hence we shall look for a substitute of this inequality in dimension $n=2$.

Due to the Parseval identity we get
$$\int \bar f (x)\partial_j g(x) \frac{x_j}{|x|} \hbox{ } dx=
\int |D|^{\frac 12} \left (  \bar f(x) \frac{x_j}{|x|}\right)
|D|^{-\frac 12}\left (\partial_j g(x)\right)\hbox{ } dx, $$
where $|D|\equiv \sqrt {-\Delta}$
and $\partial_j \equiv \frac{\partial}{\partial x_j}$,
and then due to the Cauchy--Schwartz inequality we deduce
\begin{equation}\label{kpv14}\int \bar f (x)\partial_j g(x) 
\frac{x_j}{|x|} \hbox{ } dx
\leq \left \| 
|D|^{\frac 12} \left ( \bar f(x) \frac{x_j}{|x|}\right)
\right \|_{L^2_x}
\||D|^\frac 12 g\|_{L^2_x}\end{equation}
$$=  \left \| 
|D|^{\frac 12} \left ( \bar f(x) \frac{x_j}{|x|}\right)
\right \|_{L^2_x}
\|g\|_{\dot H^\frac 12_x}.$$

On the other hand we have the following chain of inequalities:
\begin{equation}\label{basisss} \left \| 
|D|^{\frac 12} \left ( \bar f(x) \frac{x_j}{|x|}\right)
\right \|_{L^2_x}\leq 
\left \| 
|D|^{\frac 12 } \left ( \bar f(x) \frac{x_j}{|x|}\right)
- \bar f |D|^{\frac 12} \left (\frac{x_j}{|x|}\right)
\right \|_{L^2_x}\end{equation}
$$
+\left \| \bar f |D|^{\frac 12} \left (\frac{x_j}{|x|}\right)
\right \|_{L^2_x}\leq \|f\|_{\dot H^{\frac 12}_x}
+\left \| \bar f |D|^{\frac 12} \left (\frac{x_j}{|x|}\right)
\right \|_{L^2_x}$$
where we have used the inequality
$$\||D|^s (f g)- f (|D|^s g)\|_{L^2_x}\leq 
\|g\|_{L^\infty_x}\||D|^s f\|_{L^2_x}
$$
(for a proof see \cite{kepove}).

On the other hand the following inequality can be proved:
$$\left ||D|^{\frac 12} \left (\frac{x_j}{|x|}\right )\right |
\leq \frac{C}{|x|^{\frac 12}},$$
(for a proof see for example \cite{p})
and due to the Sobolev embedding it implies
$$\left \| \bar f |D|^{\frac 12 } \left (\frac{x_j}{|x|}\right)
\right \|_{L^2_x}\leq C \left \|\frac{1}{|x|^{\frac 12}} 
\bar f\right \|_{L^2_x}
\leq C \|f\|_{\dot H^{\frac 12}_x}.$$

By combining this inequality with \eqref{basisss} we get:
$$\left \| 
|D|^{\frac 12} \left ( \bar f(x) \frac{x_j}{|x|}\right)
\right \|_{L^2_x}\leq 
C  \|f\|_{\dot H^{\frac 12}_x},$$
that in turn with \eqref{kpv14} gives
$$|a(f, g)|\leq \|f\|_{\dot H^{\frac 12}_x} \|g\|_{\dot H^{\frac 12}_x}.
$$

The proof is complete.

\hfill$\Box$

\begin{lem}
Let $u(t,x)$ be the solution to 
\eqref{cau} where $n\geq 1$ and
$V(x)$ 
satisfies the same assumptions as in theorem \ref{mainUCSL0},
then the following a --priori estimates are satisfied:
\begin{equation}\label{energyhalf}
\|u(t)\|_{\dot H^\frac 12_x}^2
\leq \|u(t)\|_{\dot H^\frac 12_V}^2\leq \|f\|_{\dot H_V^\frac 12}^2 
\hbox{ } \forall t\in \mathbf R.
\end{equation}
Moreover we have
\begin{equation}\label{sdensity}\|u(t)\|_{L^2_x}^2 +\| u(t)\|_{\dot H^1_x}^2
\leq C(\|f\|_{L^2_x}^2 +\| f\|_{\dot H^1_x}^2) \hbox{ } \forall t\in \mathbf R.
\end{equation}

\end{lem}

\noindent {\bf Proof.}
Due to \eqref{charge} and
\eqref{perturbedenergy} we get:
$$\|u(t)\|_{L^2_x}= \|f\|_{L^2_x} \hbox{ and }
\|u(t)\|_{\dot H^1_V}= \|f\|_{\dot H^1_V}.$$

Hence the r.h.s. in \eqref{energyhalf} follows 
by interpolation (see \cite{tr}).

Next notice that by
hypothesis $V(x)\geq 0$ and then 
$$\|h\|_{\dot H^1_x}^2\leq \int (|\nabla h(x)|^2 + V(x)|h(x)|^2) \hbox{ } dx
=\|h\|_{\dot H^1_V}^2.$$

Hence the l.h.s. in \eqref{energyhalf} will follow 
again from an interpolation argument.

\vspace{0.1cm}

The proof of \eqref{sdensity}
follows from \eqref{charge} and \eqref{energy}.

\hfill$\Box$ 

\vspace{0.1cm}

\noindent {\bf Proof of theorem \ref{invariance}.}
Due to the r.h.s. in \eqref{energyhalf} it is easy to verify
that \eqref{famboc}
will follow from the following inequality:
\begin{equation}\label{cern}
\int |x|  |u(t)|^2 \hbox{ } dx
\leq \int |x|  |f(x)|^2 \hbox{ } dx
+ C |t| \|f\|_{\dot H^\frac 12_V}^2.\end{equation}

In the sequel we shall prove \eqref{cern} only for $t>0$ (the proof
is similar in the case $t<0$).

\vspace{0.1cm}

{\em Proof of \eqref{cern} for $t>0$}

\vspace{0.2cm}

\noindent  Since $u(t, x)$ solves \eqref{cau} 
we have the following identity: 
\begin{equation}\label{viriel}\frac d{dt} \int |x| |u(t)|^2 \hbox{ } dx
= \int |x| (\partial_t u(t) \bar u(t) 
+ u(t) \partial_t \bar u(t) )
\hbox{ } dx\end{equation}$$
= 2 {\mathcal Im}\int |x| \Delta u(t) \bar u(t) \hbox{ } dx=
-2 {\mathcal Im} \int \bar u (t) 
\nabla u(t) \cdot \frac x{|x|}\hbox{ } dx$$
where we have used integration by parts.

In order to simplify the notation we introduce the function
\begin{equation}\label{G}G(t)\equiv -2 {\mathcal Im} 
\int \bar u (t) \nabla u(t) 
\cdot \frac x{|x|}\hbox{ } dx
\end{equation}
that due to \eqref{furt} and \eqref{energyhalf}
satisfies:
\begin{equation}\label{furt2} |G(t)|\leq C
\|f\|_{\dot H^\frac 12_V}^2 \hbox{ } \forall t \in \mathbf R
\end{equation}
(notice that in fact in order to justify this computation
we have to work by density and first we  assume $f\in L^2_x \cap \dot H^1_x$.
In this way due to \eqref{sdensity} we have that
$u(t)\in L^2_x \cap \dot H^1_x$
and hence we are in position to apply \eqref{furt} with $h=u(t)$). 

By combining this inequality with \eqref{viriel} we get
\begin{equation*}\int |x||u(t)|^2 \hbox{ } dx 
\leq \int |x| |u(0)|^2 \hbox{ } dx
+ C t \|f\|_{\dot H^\frac 12_V}^2\end{equation*}
$$=\int |x| |f(x)|^2 \hbox{ } dx
+ C t  \|f\|_{\dot H^\frac 12_V}^2 \hbox{ } \forall t>0.$$

\vspace{0.1cm}

{\em Proof of \eqref{secunda}}

\vspace{0.2cm}

\noindent We shall prove \eqref{secunda} only in the case 
$t\rightarrow \infty$
(the case $t\rightarrow -\infty$ can be treated in a similar way).
Notice that \eqref{viriel} implies:
\begin{equation}\label{opul}\int \frac {|x|}{t}|u(t)|^2 \hbox{ } dx
= \int \frac {|x|}{t}|f(x)|^2 \hbox{ } dx + 
\frac{\int_0^t G(s) \hbox{ } ds}t
\hbox{ } \forall t>0.\end{equation}

On the other hand  \eqref{somcon77} and \eqref{dot7} 
(where we choose $\psi\equiv |x|$) imply
\begin{equation}\label{stam}
\lim_{t\rightarrow \pm \infty} G(t)= \pm 2 \|f\|_{\dot H^\frac 12_V}^2.
\end{equation}

By combining  \eqref{furt2}, \eqref{opul} and \eqref{stam} we finally get
\begin{equation}\label{massi}\lim_{t\rightarrow \infty}
\int \frac {|x|}{t}|u(t)|^2 \hbox{ } dx= 2\|f\|_{\dot H^\frac 12_V}^2.
\end{equation}

\hfill$\Box$

\section{Proof of theorems \ref{mainSL} and \ref{mainUCSL}}
\label{benfra}
We split the proof in two steps.

\vspace{0.1cm}

{\em Proof of r.h.s. in \eqref{eqSL}}

\vspace{0.2cm}

\noindent It is sufficient to consider the identity \eqref{vai}
where the generic function $\psi$ is replaced
with the family of rescaled functions $R \phi\left (\frac xR \right)$ and 
$\phi(x)\equiv \langle x \rangle$.
In fact notice that the function $\langle x \rangle$ 
is convex, increasing and moreover
$$-\Delta^2 (\langle x \rangle)\geq 0 \hbox{ } 
\forall x\in {\mathbf R}^n, n\geq 3$$
as a direct computation shows.

Let us point--out that the l.h.s.
in \eqref{eqSL}
follows from theorem \ref{mainUCSL}.

\vspace{0.1cm}

{\em Proof of theorem \ref{mainUCSL}} 

\vspace{0.2cm}

\noindent Next we shall make use of
the following identity 
\begin{equation}\label{cube}\nabla \bar u D^2 \psi \nabla u=
\partial_{|x|}^2\psi |\partial_{|x|} u|^2 + \frac{\partial_{|x|} 
\psi}{|x|}|\nabla_\tau u|^2,\end{equation}
where $\psi$ is a radially symmetric function
and $u$ is a generic function.

\noindent First of all let us notice that if we choose in the identity \eqref{vai}
the function $\psi\equiv |x|$ then we get:
\begin{equation*}
\int \int_{|x|>1} \frac{|\nabla_\tau u|^2}{|x|} \hbox{ } dxdt <\infty,
\end{equation*}
where we have used
\eqref{ngeq4} and \eqref{cube}.

In particular we get
\begin{equation}\label{revised}
\lim_{R\rightarrow \infty}
\int \int_{|x|>R} \frac{|\nabla_\tau u|^2}{|x|} \hbox{ } dx=0.
\end{equation}

For any $k\in \mathbf N$
we fix a function
$h_k(r)\in C^\infty_0(\mathbf R; [0, 1])$ such that:
\begin{equation}\label{hk}h_k(r)=1 \hbox{  } \forall r\in {\mathbf R}
\hbox{ s.t. } |r|<1, h_k(r)=0 \hbox{  } \forall r \in {\mathbf R}
\hbox{ s.t. } |r|>\frac{k+1}{k}, 
\end{equation}
$$h_k(r)= h_k(-r) \hbox{  } \forall r \in \mathbf R.$$

Let us introduce the functions
$\psi_k(r), H_k(r) \in C^\infty( {\mathbf R})$:
\begin{equation}\label{kr}\psi_k(r) = \int_0^r (r-s) h_k(s) ds \hbox{  } \hbox{ and } \hbox{  }
H_k(r)=\int_0^r h_k(s) ds.\end{equation}

Notice that 
\begin{equation}\label{differ}
\psi_k''(r)=h_k(r), \psi_k'(r)= H_k(r) \forall r\in {\mathbf R} \hbox{ and }
\lim_{r\rightarrow \infty} \partial_r \psi_k(r)=\int_0^\infty h_k(s)ds.
\end{equation}

Moreover an elementary computation 
shows that: 

\begin{equation}\label{bilaplac}\Delta^2 \psi_k(x)= 
\frac C{|x|^3} \hbox{  } \forall x\in 
{\mathbf R}^n \hbox{ s.t. } |x| \geq 2 \hbox{ and } n\geq 4,
\end{equation}
where $\Delta^2$ is the bilaplacian operator.

Thus the functions
$\phi=\psi_k$ satisfy the assumptions of propositions \ref{noma}
and \ref{noma3}.

\vspace{0.1cm}

In the sequel we shall need the rescaled functions
\begin{equation}\label{kR}\psi_{k, R}(x)=R \psi_k\left (\frac{x}{R} \right ) \forall x\in {\mathbf R}^n, k\in 
{\mathbf N} \hbox{ and }R>0,
\end{equation}
where $\psi_k$ is defined in \eqref{kr}.
Notice that by combining the general identity 
\eqref{cube} with \eqref{vai}, where we choose $\psi= 
\psi_{k,R}$ defined in \eqref{kR}, and recalling \eqref{differ} we get:
\begin{equation}\label{limit}\int \int
\left [ \partial_{|x|}^2 \psi_{k,R} |\partial_{|x|} u|^2 + \frac{\partial_{|x|} \psi_{k,R}}
{|x|} | \nabla_\tau u|^2
\right .\end{equation}
\begin{equation*}
\left . -\left (\frac 14
\Delta^2 \psi_{k,R} + \partial_{|x|}V \partial_{|x|} \psi_{k,R}\right )
|u|^2 \right ]dx dt 
= \left (\int_0^\infty h_{k}(s) ds
\right ) \|f \|_{\dot H^\frac 12_V }^2 \forall k\in {\mathbf N}, R>0.
\end{equation*}

By using \eqref{uno}, \eqref{due3}
and \eqref{revised} we get:
\begin{equation*}
\lim_{R\rightarrow \infty} \int\int
\left [\partial_{|x|} \psi_{k,R}\frac{|\nabla_\tau u|^2}{|x|} 
- \left (\frac 14 \Delta^2 \psi_{k,R}  + \partial_{|x|}V \partial_{|x|} \psi_{k,R}\right )
|u|^2 \right ]dx dt=0
\end{equation*}
for every $k\in {\mathbf N}$.

We can combine this fact with \eqref{limit} in order to deduce:

\begin{equation}\label{crusia}\lim_{R\rightarrow \infty} \int\int
\partial_{|x|}^2\psi_{k,R} |\partial_{|x|} u|^2
dx dt =\left (\int_0^\infty h_{k}(s) ds
\right ) \|f \|_{\dot H^\frac 12_V }^2 \forall k\in {\mathbf N}.
\end{equation}

On the other hand,
due to the properties of $h_{k}$ (see \eqref{hk}), we get
$$\frac 1R \int \int_{B_R}|\partial_{|x|} u|^2
dx dt \leq \int \int
\partial_{|x|}^2 \psi_{k,R}  |\partial_{|x|} u|^2  dt dx $$$$
= \frac 1R \int \int
h_k\left( \frac {x}R\right)|\partial_{|x|} u|^2  dt dx 
\leq \frac 1R \int\int_{|x|<\frac {k+1}{k}R} |\partial_{|x|} u|^2 dxdt$$
that due to \eqref{crusia} implies:
\begin{equation}\label{limsup}\limsup_{R\rightarrow \infty} \frac 1R \int 
\int_{|x|<R}|\partial_{|x|} u|^2 dxdt 
\leq  \left (\int_0^\infty h_k(s) ds \right ) \|f\|_{\dot H^\frac 12_V}^2
\end{equation}$$\leq\frac {k+1}{k}
\liminf_{R\rightarrow \infty} \frac 1R \int\int_{|x|<R}|\partial_{|x|} u|^2 dxdt
\hbox{  } \forall k \in {\mathbf N}.$$

Since $k\in \mathbf N$ is arbitrary and since 
the following identity is trivially satisfied:
$$\lim_{k\rightarrow \infty}\int_0^\infty h_k(s) ds=1,$$
we can deduce easily
\eqref{eqUCSL} by using \eqref{limsup}.

\vspace{0.1cm}

The proof is complete.

\hfill$\Box$

\section{Proof of theorems \ref{mainSL3} and \ref{mainUCSL3}}
\label{fanucci}
The proofs are similar in principle to the ones of theorems \ref{mainSL}
and \ref{mainUCSL}, except that in dimension $n=3$
we cannot use \eqref{bilaplacian} and hence lemma \ref{noma}, that have been
proved only in dimension $n\geq 4$.
However for completeness we shall give in the sequel the details
of the proof in dimension $n=3$.
\\
\\
We split the proof in two steps.

\vspace{0.1cm}

{\em Proof of r.h.s. in \eqref{eqSL3}}

\vspace{0.2cm}

\noindent It is sufficient to consider the identity \eqref{vai}
where the generic function $\psi$ is replaced
with the family of rescaled functions $R \phi\left (\frac xR \right)$ and 
$\phi(x)\equiv \langle x \rangle$.
In fact notice that the function $\langle x \rangle$ 
is convex, increasing and moreover
$$-\Delta^2 (\langle x \rangle)\geq 0 \hbox{ } 
\forall x\in {\mathbf R}^3$$
as a direct computation shows.

Let us point--out that the l.h.s.
in \eqref{eqSL3}
follows from theorem \ref{mainUCSL3}.

\vspace{0.1cm}

{\em Proof of theorem \ref{mainUCSL3}} 

\vspace{0.2cm}

\noindent Following the proof of \eqref{revised} we get:
\begin{equation}\label{revised3}
\lim_{R\rightarrow \infty}\int \int_{|x|>R} \frac{|\nabla_\tau u|^2}{|x|} \hbox{ } dxdt=0,
\end{equation}
(in this case of course we have to use \eqref{ngeq3}
instead of \eqref{ngeq4}).

We fix a function
$h(r)\in C^\infty_0(\mathbf R; [0, 1])$ such that:
$$h(r)=1 \hbox{  } \forall r\in {\mathbf R}
\hbox{ s.t. } |r|<\frac 12, h(r)=0 \hbox{  } \forall r \in {\mathbf R}
\hbox{ s.t. } |r|>1, 
$$
$$h(r)= h(-r) \hbox{  } \forall r \in \mathbf R.$$

Let us introduce the functions
$\psi(r), H(r) \in C^\infty( {\mathbf R})$:
\begin{equation}\label{benc}
\psi(r) = \int_0^r (r-s) h(s) ds \hbox{  } \hbox{ and } \hbox{  }
H(r)=\int_0^r h(s) ds.
\end{equation}

Notice that 
\begin{equation}\label{differ3}
\psi''(r)=h(r), \psi'(r)= H(r) \forall r\in {\mathbf R} \hbox{ and }
\lim_{r\rightarrow \infty} \partial_r \psi(r)=\int_0^\infty h(s)ds.
\end{equation}

Moreover an elementary computation 
shows that: 

\begin{equation}\label{bilaplac3}\Delta^2 \psi(x)= 0 \hbox{  } \forall x\in 
{\mathbf R}^3 \hbox{ s.t. } |x| \geq 1 ,
\end{equation}
where $\Delta^2$ is the bilaplacian operator (recall that we are working in dimension
$n=3$) and $\psi$ is defined in \eqref{benc}

Notice also that the function
$\psi$ given above satisfies the assumptions of proposition \ref{noma3}.
In the sequel we shall need the rescaled functions
\begin{equation}\label{fort}\psi_{R}(x)=R \psi\left (\frac{x}{R} \right ) \forall x\in {\mathbf R}^3
\hbox{ and }R>0.\end{equation}

By combining the identity \eqref{cube} with \eqref{vai}, where we choose $\psi= 
\psi_{R}$, and recalling \eqref{differ3} we get:
\begin{equation}\label{limit3}\int \int
\left [ \partial_{|x|}^2 \psi_{R} |\partial_{|x|} u|^2 + \frac{\partial_{|x|} \psi_{R}}
{|x|} | \nabla_\tau u|^2
\right .\end{equation}
\begin{equation*}
\left . -\left (\frac 14
\Delta^2 \psi_{R} + \partial_{|x|}V \partial_{|x|} \psi_{R}\right )
|u|^2 \right ]dx dt 
= \left (\int_0^\infty h(s) ds
\right ) \|f \|_{\dot H^\frac 12_V }^2 \forall R>0.
\end{equation*}

By using \eqref{due3}
and \eqref{revised3} we get:
\begin{equation*}
\lim_{R\rightarrow \infty} \int \int
\left (\partial_{|x|} \psi_{R}\frac{|\nabla_\tau u|^2}{|x|} 
- \partial_{|x|}V \partial_{|x|} \psi_{R}
|u|^2 \right )dx dt=0.
\end{equation*}

We can combine this fact with \eqref{limit3} in order to deduce:

\begin{equation}\label{crusia3}\lim_{R\rightarrow \infty} \int\int
\left ( \partial_{|x|}^2\psi_{R} |\partial_{|x|} u|^2
 - \frac 14
\Delta^2 \psi_{R}  |u|^2 \right )dx dt =\left (\int_0^\infty h(s) ds
\right ) \|f \|_{\dot H^\frac 12_V }^2.
\end{equation}

On the other hand,
due to the definition of the functions $h(|x|)$ and $\psi_R(|x|)$ (see \eqref{benc}
and \eqref{fort}), we get
$$\int \int
\left( \partial_{|x|}^2 \psi_{R}  |\partial_{|x|} u|^2  
- \frac 14
\Delta^2 \psi_{R}  |u|^2 \right ) 
dx dt $$$$
\leq  \int \int \left (\frac 1R 
h\left( \frac {x}R\right)|\partial_{|x|} u|^2 
- \frac 1{4 R^3} 
(\Delta^2 \psi)^- \left( \frac xR \right)  |u|^2 \right ) dx dt
$$
(here $(\Delta^2 \psi)^-$ represents the negative part of $\Delta^2 \psi$),
that in turn due to \eqref{bilaplac3} and again to the property 
of the support of $h(|x|)$ implies:
$$\int \int
\left( \partial_{|x|}^2 \psi_{R}  |\partial_{|x|} u|^2  
- \frac 14
\Delta^2 \psi_{R}  |u|^2 \right) 
dx dt $$
$$
\leq \frac 1R  \int \int_{|x|<R}
\left ( |\partial_{|x|} u|^2 
+ \frac C {R^2} 
|u|^2 \right )dx dt.
$$

Finally we can combine this estimate 
with \eqref{crusia3} in order to deduce
\eqref{eqUCSL3}.

\vspace{0.1cm}

The proof is complete.

\hfill$\Box$
 
\section{Appendix}\label{appendix}

In this section we shall
work with potentials $V(x)\in L^\infty_x$ such that
\begin{equation}\label{rageweak}
\lim_{|x|\rightarrow \infty} V(x)=0 \hbox{ and }
\lim_{|x|\rightarrow \infty} |x||\partial_{|x|} V(x)|=0. 
\end{equation}

We point out that 
this appendix does not contain any essential novelty,
since the content of proposition \ref{ragetheorem} in part follows
from the results in \cite{mourre} (see remark \ref{remmourre}).

However, our aim in this appendix is to show in a very simple way how to deduce 
a stronger version
of the usual RAGE theorem by using as a starting point the classical RAGE theorem himself
(see proposition \ref{ragetheorem}) and
by avoiding the use of the general Mourre theorem in \cite{mourre}.

\vspace{0.1cm}

It is well--known that to every bounded potential $V(x)$,
we can associate a corresponding splitting of the Hilbert space $L^2_x$
as the direct sum
of the projection onto the continuous spectrum and onto the pure point spectrum,
that shall be denoted respectively as 
$L^2_c$ and $L^2_{pp}.$

The space $L^2_c$  can be splitted in turn as the direct sum
of the projection onto the singular spectrum  
and onto the absolutely continuous spectrum,
that shall be denoted respectively as 
$L^2_s$ and $L^2_{ac}.$

\vspace{0.1cm}

Hereafter we shall make use of the following version of the 
RAGE theorem:
\begin{equation}\label{ragechi}
\hbox{ if } u\in {\mathcal C}_t (L^2) \hbox{ solves \eqref{cau} where } 
f \in L^2_c
\end{equation}
\begin{equation*}
\hbox{ then } \lim_{T\rightarrow \infty}
\frac 1T \int_{-T}^T\int_{|x|<R} |u(t)|^2 \hbox{ } dxdt=0
\hbox{ } \forall R>0.
\end{equation*}

Actually the classical RAGE theorem is much more general 
than \eqref{ragechi} (see \cite{RS2}),
however \eqref{ragechi} will be enough for our purposes.

\vspace{0.1cm}

One of the aims of this appendix is to show
how the RAGE theorem
implies \eqref{or}, under the decay assumptions 
on $V(x)$ given in \eqref{rageweak}.
Let us recall that \eqref{or} is a stronger
version of a result obtained in \cite{Enss}
(see remark \ref{Enssremark}).

Another point in this appendix is to show how it is possible to 
prove a decay of the solution $u$ pointwisely in time,
by using as a starting point 
the decay given in \eqref{ragechi}.

Next we state the main result of the section.
\begin{prop}\label{ragetheorem}
Let $V(x)$ be a function that satisfies \eqref{rageweak}.
Assume that the point spectrum of $-\Delta+V$ is 
the empty set.
Then we have:
\begin{equation}\label{or2}
\lim_{t\rightarrow \pm \infty}
 \left \|\frac xt u(t) - 2{\bf i} \nabla u(t)
\right \|_{L^2_x}=0
\end{equation}
where
$$u(t)\equiv e^{{\bf i}t\Delta_V}f
\hbox{ with } \int |x|^2 |f|^2 \hbox{ } dx$$
and
\begin{equation}\label{mainappendix}\lim_{t\rightarrow \pm \infty}
\int W(x) |u(t)|^2 \hbox{ } dx=0
\end{equation}
where
$$u(t)\equiv e^{{\bf i}t\Delta_V}f
\hbox{ with } f\in L^2_x$$
and
$W(x)\in L^\infty_x$ satisfies
$$\lim_{|x|\rightarrow \infty} W(x)=0.$$
In particular we have:
\begin{equation}\label{weaklimit}e^{{\bf i}t\Delta_V}f\rightharpoonup
0
 \hbox{ as } t\rightarrow \pm \infty \hbox{ } \forall f\in L^2_x.
\end{equation}
\end{prop} 

\begin{remark}\label{remmourre}
Looking at the proof of the RAGE theorem, which is based on
Wiener's result about the decay of the Fourier transform 
of a measure, it is no difficult to show that
\eqref{mainappendix} is trivially satisfied
provided that the projection on the absolutely continuous spectrum 
$L^2_{ac}$ coincides with $L^2$, or equivalently 
when the projection on the singular spectrum $L^2_s$ is trivial. 
Actually this fact has been proved in \cite{mourre}, provided that
the potential $V(x)$ 
satisfies \eqref{rageweak}.
However we have decided to present our own proof
of proposition \ref{ragetheorem} due to its simplicity.
\end{remark}

\noindent{\bf Proof of proposition \ref{ragetheorem}.}
For simplicity we shall treat only the limit as $t\rightarrow \infty$ 
(the case $t\rightarrow -\infty$ can be studied  in a similar way).
 
\vspace{0.1cm}

{\em Proof of \eqref{or2}}

\vspace{0.2cm}
\noindent We have that
$L^2_x \equiv L^2_{c}.$
Hence we can use the 
RAGE theorem (see \eqref{ragechi}) in order to deduce
that 
\begin{equation}\label{pallar}
\lim_{T\rightarrow \infty} \frac 1T \int_{-T}^T\int_{|x|<R} |u(t)|^2 \hbox{ } dxdt=0,
\end{equation}
where
$u(t)=e^{{\bf i}t\Delta_V}f$.

Notice that $u\equiv e^{{\bf i}t\Delta_V}f\in 
{\mathcal C}_t(L^2_x)$ can be also defined as the unique solution to
$$
{\bf i}\partial_t u - \Delta u + V(x) u=0,$$
$$u(0)=f,$$
hence we can use the estimate \eqref{ragecaze}
given in lemma \ref{ragecazenave} in order to get:
$$\left \|\frac xt u(t) - 2{\bf i} \nabla u(t)
\right \|_{L^2_x}^2$$
$$\leq \frac{\int |x|^2 |f(x)|^2  \hbox{ } dx}{t^2}
+\frac{8 \int_0^t \left(\int U(x) |u(s)|^2\hbox{ } dx \right ) \hbox{ } ds}{t}$$
where
$$U(x)\equiv 8\left(V(x) + \frac 12 |x| \partial_{|x|}V(x)\right ).$$

Due to \eqref{pallar} and to the decay  assumption \eqref{rageweak}, 
the previous estimate 
can be rewritten as follows:
\begin{equation}\label{metaf5}\left \|\frac xt u(t) - 2{\bf i} \nabla u(t)
\right \|_{L^2_x}^2\end{equation}
$$\leq \frac{\int |x|^2 |f(x)|^2  \hbox{ } dx}{t^2}
+\varphi(t)$$
where
$$\lim_{t\rightarrow \infty} \varphi(t)=0.$$

Hence the proof of \eqref{or2} follows easily.

\vspace{0.2cm}
{\em Proof of \eqref{mainappendix}} 
\vspace{0.1cm}

\noindent Due to \eqref{charge}
and to a density argument, 
it is sufficient to prove \eqref{mainappendix} under the assumption
$f\in C^\infty_0({\mathbf R}^n)$.

Next notice that \eqref{metaf5} can be written as 
\begin{equation}\label{metaf}4 \left \|\nabla (e^{{\bf i}\frac{|x|^2}{4t}} u(t))
\right \|_{L^2_x}^2\leq \frac{\int |x|^2 |f(x)|^2  \hbox{ } dx}{t^2}
+\varphi(t)\end{equation}
where
$$\lim_{t\rightarrow \infty} \varphi(t)=0.$$

Since we are assuming $f\in C^\infty_0({\mathbf R}^n)$
we have that the R.H.S. in
\eqref{metaf} goes to zero as $t\rightarrow \infty$.

Assume $n>2$. By the Sobolev embedding
$$\dot H^1 ({\mathbf R}^n)\subset L^{2^*}({\mathbf R}^n)
\hbox{ where } 2^*\equiv \frac{2n}{n-2},$$
we get
\begin{equation}\label{sobolevski}
\lim_{t\rightarrow \infty} \|u(t)\|_{L^{2^*}}=0.
\end{equation}

This estimate
in conjunction with the H\"older inequality implies:
$$\int_{|x|<R} |u(t)|^2 \hbox{ } dx
\leq C R^2 \|u(t)\|_{L^{2^*}}^2  \hbox{ } \forall R>0,$$
that in turn can be combined with \eqref{sobolevski}  in order to give
\begin{equation}\label{l2localized}\lim_{t\rightarrow \infty}
\int_{|x|<R} |u(t)|^2 \hbox{ } dx=0 \hbox{ } \forall R>0.
\end{equation}

Notice that \eqref{mainappendix} 
follows by 
\eqref{l2localized}, \eqref{charge} and the decay assumption 
of $W(x)$ at infinity.

Finally notice that \eqref{weaklimit}
follows by combining \eqref{charge}
with \eqref{mainappendix}.

For $n=1$ and $n=2$ a similar argument works  using Gagliardo--Nirenberg inequalities instead of
the Sobolev embedding.

\hfill$\Box$

The aim of the next proposition is to show that, despite
to \eqref{mainappendix}, in general there is 
a--priori  no rate of decay for  the $L^2$ localized norm of the solutions
$u$ to \eqref{cau}.

\begin{prop}\label{noratedecay}
Assume that $V(x)\in L^\infty_x$
is any bounded potential
(possibly $V(x)\equiv 0$). Let $R>0$ be a fixed positive number
and $\gamma\in {\mathcal C} ([0, \infty); {\mathbf R})$
be any function such that
$$\lim_{t\rightarrow \infty}\gamma (t)=\infty.$$
Then there exists $g\in L^2_x$ 
(that depends on $R, V(x)$ and $\gamma(t)$)such that
$$\int_{|x|<R}|u(t_n)|^2 \hbox{ } dx> \frac n{\gamma(t_n)}$$
where $\{t_n\}_{n\in \mathbf N}$ is a suitable sequence 
$\lim_{n\rightarrow \infty} t_n=\infty$
and
$u(t)\in {\mathcal C}_t(L^2_x)$
satisfies
\begin{equation}\label{CPappendix}
{\bf i}\partial_t u - \Delta u + V(x) u=0,
\end{equation} $$u(0)=g.$$
\end{prop}

\noindent{\bf Proof.}
We claim the following fact:

\begin{equation}\label{claim} 
\|{\mathcal U}_{R,V}(t)\|_{{\mathcal L}(L^2_x, L^2_x)}
\equiv 1 \
\hbox{ } \forall t\in {\mathbf R} \hbox{ where }
\end{equation}$${\mathcal U}_{R,V}(t): L^2_x\ni g \rightarrow \chi_{\{
|x|<R\}}u(t)\in L^{2}_x
$$$$\hbox{ and $u(t)$ denotes the unique solution to the Cauchy problem
\eqref{CPappendix}}.
$$

Notice that due to \eqref{claim} we get:
$$\lim_{t\rightarrow \infty}
\gamma(t)\|{\mathcal U}_{R,V}(t)\|_{{\mathcal L}(L^2_x, L^{2}_x)}=\infty$$
and in particular
due to the Banach-Steinhaus theorem the operators
$\gamma(t) {\mathcal U}_{R,V}(t)$ cannot be pointwisely bounded
or in an equivalent way there exists at least one $g\in L^2_x$ such that
\begin{equation}\label{BS}\sup_{[0, \infty)} \gamma(t)\| 
{\mathcal U}_{R,V}(t) g\|_{L^{2}_x}=\infty.
\end{equation}

On the other hand the function $t\rightarrow \gamma(t)\|{\mathcal U}_{R,V}(t)g\|_{L^{2}_x}$
is bounded on bounded sets of $[0, \infty)$
and hence \eqref{BS} implies that
$$\limsup_{t\rightarrow \infty} \gamma(t)\|{\mathcal U}_{R,V}(t)
g\|_{L^{2}_x}=\infty$$
which completes the proof.

\vspace{0.1cm}

Next we shall prove \eqref{claim}.
Let us fix any function $f_R\in L^2_x$
such that $$\|f_R\|_{L^2_x}\equiv 1 \hbox{ and } supp \hbox{ }
f_R\subset \{|x|<R\}.$$

Notice that we have
\begin{equation}\label{trunc}
{\mathcal U}_{R,V}(t)e^{-{\bf i}t\Delta_V} f_R \equiv 
\chi_{\{|x|<R\}}e^{{\bf i}t\Delta_V}e^{-{\bf i}t\Delta_V} f_R\equiv
f_R\hbox{ }
\forall t\in \mathbf R
\end{equation}
where we have used the group property
of $e^{{\bf i}t\Delta_V}$
and the assumption done on the support of $f_R$.

In particular due to \eqref{charge} and \eqref{trunc} 
we get:
$$\|{\mathcal U}_{R,V}(t)\|_{{\mathcal L}(L^2_x, L^2_x)}
\geq \|{\mathcal U}_{R,V}(t)(e^{-{\bf i}t\Delta_V} f_R)\|_{L^2_x}
\equiv \|f_R\|_{L^2_x}\equiv 1.$$

On the other hand \eqref{charge} implies 
trivially the opposite
inequality
$$\|{\mathcal U}_{R,V}(t)\|_{{\mathcal L}(L^2_x, L^2_x)}
\leq 1$$ and hence \eqref{claim} is proved. 

\hfill$\Box$

\end{document}